\documentclass{ccgrt} 

\newtheorem{thm}{Theorem}[section]
\newtheorem{cor}[thm]{Corollary}
\newtheorem{lem}[thm]{Lemma}

\numberwithin{equation}{subsection}

\newtheorem{conj}[thm]{Conjecture}
\newtheorem{quest}[thm]{Question}

\usepackage{pstricks}

\begin{document}

\title{Some recent approaches in 4--dimensional surgery theory\,\footnote{\,
2000 \uppercase{M}ath. \uppercase{S}ubj. \uppercase{C}lassif. 
\uppercase{P}rimary: \uppercase{57R67}, \uppercase{57P99} 
\uppercase{S}econdary: \uppercase{55N45}. 
\newline {\it \uppercase{K}eywords:\,} 
\uppercase{P}oincar\'{e} complex, chain duality, controlled topology, 
controlled surgery, 
homology manifold.}}

\author{Friedrich Hegenbarth 
}

\address{Department of Mathematics,
University of Milano, Via C. Saldini 50, Milano, Italy 02130.
\vspace{1mm}\\
E-mail:  friedrich.hegenbarth@mat.unimi.it }

\author{Du\v san Repov\v s\,\footnote{\,\uppercase{W}ork supported by the 
\uppercase{S}lovenian
\uppercase{R}esearch
\uppercase{A}gency
grant\, \uppercase{N}o. 
{\uppercase{P}1-0292-0101-04.}} }

\address{Institute of Mathematics, Physics and Mechanics,
University of Ljubljana, P. O. Box 2964, Ljubljana, Slovenia 1001.
\vspace{1mm}\\
E-mail: dusan.repovs@fmf.uni-lj.si }

\maketitle

\abstracts
{It is well-known that an n-dimensional Poincar\'{e} complex $X^n$, 
$n \ge 5$, has the homotopy type of a compact topological 
$n$-manifold if 
the total surgery obstruction $s(X^n)$ vanishes. 
The present paper discusses recent attempts to prove analogous result in 
dimension 4. We begin by reviewing the necessary
algebraic and controlled surgery theory. 
Next, we discuss the key idea of Quinn's approach. Finally, we present 
some cases of special fundamental 
groups, due to the authors and to Yamasaki.}

\setcounter{section}{0}


\section{Introduction} 

Classical surgery methods of 
Browder--Novikov--Wall break down in dimension $4$. 
The
Wall groups, depending only on the fundamental
group, do not seem to be strong enough as obstruction
groups to completing the surgery. It is strongly 
believed that for free nonabelian fundamental
groups of rank $r \ge 2$ Wall groups are not
sufficient
(\cite{FrQu}). 
Nevertheless, one can make progress using
controlled surgery theory to produce controlled
embeddings of $2$--spheres needed for surgery, 
by using results of Quinn (\cite{Qu3}). However, 
this works only if the control map satisfies
the $UV^1$--condition (in fact, one needs the $UV^1(\delta)$--condition, for sufficiently
small $\delta > 0$). 
The obstructions belong to a
controlled Wall group. Its construction is
conceptual, done by means of 
Ranicki's machinery
(\cite{PQR})
(\cite{PedYa})
(\cite{Ra3})--(\cite{Ra6})
(\cite{RaYa3})
(\cite{RaYa2})
(\cite{Ya}).
 
The most important fact is that the
controlled groups are homology groups,
so they can be calculated.
The following is a basic result in
surgery theory in dimension $n \ge 5$.
It is due to Ranicki (\cite{Ra3}), see $\S2$ below:

\begin{thm}
Suppose that
$X^n$ is a Poincar\'e 
$n$--complex, $n \ge 5$,
with total surgery obstruction $s(X^n) = 0$.
Then $X^n$ is (simple) homotopy equivalent
to a topological $n$--manifold.
\end{thm}

At his talk in 2004
Quinn (\cite{Qu5}) proposed 
a strategy
to extend
Theorem 1.1 to 
dimension
$n = 4$. 
More precisely, he proposed how to prove
the following
 conjecture:

\begin{conj}
If $X^4$ is a Poincar\'e $4$--complex with $s(X^4) = 0$,
then $X^4$ is (simple) homotopy equivalent
to an ANR  homology $4$--manifold.
\end{conj}

We shall 
outline the idea of Quinn's approach in $\S4$, after having prepared
the necessary preliminaries.
We shall
also prove a "stable"
version of Conjecture 1.2:

\begin{thm}
If $X^4$ is a Poincar\'e $4$--complex with $s(X^4) = 0$,
then $X^4
\#$
$ (\mathop\#\limits_1^r S^2 \times S^2)$ is (simple)
homotopy equivalent to a topological
$4$--manifold.
\end{thm}

In the rest of $\S4$ we shall
prove special cases
of Conjecture 1.2,  when $\pi_1(X^4)$ is free nonabelian.
We begin by establishing some notations and results which are 
needed for our presentation. We acknowledge the referee for comments and suggestions.

\section{Notations and basic results of algebraic surgery}
   
Let $\Lambda$
be some ring with anti--involution.
Here we will only consider
$\Lambda = \mathbb{Z} [\pi_1]$,
the integral group ring of a fundamental group of a space with trivial orientation character.
An $n$--quadratic chain complex is a pair
$(C_{\mathop\#}, \psi_{\mathop\#})$,
where
$C_{\mathop\#}$
is a free
$\Lambda$--module chain complex and
$\psi_{\mathop\#}=\{\psi_{s} \mid s=0,1,...\}$
is a collection of
$\Lambda$--homomorphisms
$\psi_{s}: C^{n-r-s} \rightarrow C_{r}$
satisfying certain relations.
Here 
$C^{\mathop\#}$
denotes the
$\Lambda$--dual cochain complex of 
$C$
(as in (\cite{Wa})).
The pair is a quadratic $n$--Poincar\'{e} complex if
$\psi_{0} + \psi_{0}^{\mathop\#}: C^{n-{\mathop\#}} \rightarrow C_{\mathop\#}$
is a chain equivalence, where
$\psi_{0}^{\mathop\#}$
is a the dual map of
$\psi_{0}$.

There is the notion of $(n+1)$-quadratic (Poincar\'{e}) pairs,
hence of "cobordism" between $n$--quadratic Poincar\'{e} complexes,
which is an equivalence relation.
Let $L_{n}(\Lambda)$
be the set of
equivalence classes
of $n$--quadratic Poincar\'{e} complexes.
It has a group structure induced by direct sum constructions.

If $n=2k$ then particular examples of quadratic Poincar\'{e} complexes are given by surgery kernels
$[K_{k}(f),\Lambda,\mu]$
of degree 1 normal maps
$(f,b):M^{n} \rightarrow X^{n}$
(see 
(\cite{Ra1})
(\cite{Ra2})
(\cite{Ra4})
(\cite{Ra7})
).
This gives an isomorphism of
$L_{2k}(\Lambda)$
with Wall groups
$L_{2k}(\pi_{1})$.
There is also
an isomorphism in odd dimensions
in terms of formations.
There is an 
$\Omega$--spectrum
$\mathbb{L}$
with
$\pi_{n}(\mathbb{L})=L_{n}(\{1\})$
and
$\mathbb{L}_{0}=G/TOP \times \mathbb{Z}$
(see(\cite{Qu1})).
The groups
$L_{n}(\{ 1 \})$
were calculated in
(\cite{KeMi}).
We denote by
$\dot{\mathbb{L}} \rightarrow \mathbb{L}$
the connected covering spectrum, so
$\pi_{n}(\dot{\mathbb{L}})=\pi_{n}(G/TOP)$.

If $K$ is a simplicial complex, elements
$\xi \in H_{n}(K,\dot{\mathbb{L}})$
can be represented by equivalence classes of 
compatible collections of $n$--quadratic Poincar\'{e} complexes
$\{(C_{\mathop\#}(\sigma),\psi_{\mathop\#}(\sigma)) \mid \sigma \in K \}$.
Gluing these individual quadratic complexes together
gives a "global" $n$--quadratic Poincar\'{e} complex
$(C_{\mathop\#}, \psi_{\mathop\#})$
hence an element in
$L_{n}(\Lambda)$.
There results a homomorphism
$A:H_{n}(K,\dot{\mathbb{L}})\rightarrow L_{n}(\pi_{1}(K))$,
called the assembly map.
To define the structure set
${\mathbb{S}}_{n+1}(K)$,
one considers compatible collections
$\{(C_{\mathop\#}(\sigma),\psi_{\mathop\#}(\sigma)) \mid \sigma \in K \}$
of $n$--quadratic Poincar\'{e} complexes
which are "globally" contractible,
i.e. the mapping cone of the map
$\psi_{0}+\psi_{0}^{\mathop\#}:C^{n-\mathop\#}\rightarrow C_{\mathop\#}$
is contractible (see (\cite{Ra5})).

Cobordism classes of such objects build the set
${\mathbb{S}}_{n+1}(K)$.
It has a group structure coming from direct sum constructions.

\begin{thm}
The assembly map fits into the exact sequence 
$$
\cdots \to L_{n+1}(\pi_1(K)) \mathop\to\limits^{\partial} \mathbb{S}_{n+1}(K) \to H_n (K, \dot{\mathbb{L}}) 
\mathop\to\limits^A
$$
$$
\mathop\to\limits^A L_n(\pi_1(K)) \mathop\to\limits^{\partial} \mathbb{S}_n(K) \to \cdots \;\; .
$$
\end{thm}

For any geometric Poincar\'{e} duality complex
$X^{n}$
of dimension $n$, Ranicki defined its total surgery obstruction
$s(X^n) \in \mathbb{S}_n(X^n)$
(see (\cite{Ra3})).

\begin{thm}
If $n \ge 5$ then $s(X^n)=0$ if and only if $X^n$ is (simple)
homotopy equivalent to a topological $n$--manifold.
\end{thm}

The element $s(X)$
can be decribed as follows:
Suppose that $X$
is triangulated.
The fundamental chain
$[X] \in C_{n}(X)$
defines a simple chain equivalence
$$-\cap [X] \colon C^{n-\mathop\#}(\tilde{X}) \to C_{\mathop\#}(\tilde{X}), 
\ \ \  \tilde{X} \to X$$
is the universal covering,
i.e. the desuspension of the algebraic mapping cone
$S^{-1}C(-\cap [X])$ 
of
$-\cap [X]$
is contractible.

Let $\sigma^{*}$ be the dual cell of the simplex
$\sigma \in X$
with respect of its barycentric subdivision.
The global fundamental cycle 
$[X]$
defines local cycles
$[X(\sigma)] \in C_{n-|\sigma|}(\sigma^{*}, \partial{\sigma^{*}})$,
hence it maps
$$-\cap [X(\sigma)] \colon C^{n-r-|\sigma|}(\sigma^{*}) 
\to C_{r}(\sigma^{*},\partial{\sigma^{*}}).$$
The collection
$$\{ D_{\mathop\#}(\sigma) =
S^{-1}C(-\cap [X(\sigma)])
\mid \sigma \in X \}$$
assembles to
$D=S^{-1}C(-\cap [X])$.
There are $(n-1-|\sigma|)$--quadratic Poincar\'{e}
structures
$\psi_{\mathop\#}(\sigma)$
on
$D(\sigma)$,
giving rise to an 
$(n-1)$--quadratic structure on $D$.
Then $s(X)$
is represented by the class of the compatible collection
$\{(D_{\mathop\#}(\sigma),\psi_{\mathop\#}(\sigma)) \mid \sigma \in X \}$.

\section{Controlled $L$--groups of geometric
quadratic complexes and the controlled
surgery sequence}

Geometric modules were introduced by Quinn
(\cite{Qu2})
(\cite{Qu0})
(\cite{Qu4})
(\cite{Qu6})
(see also (\cite{RaYa1})). 
We introduce here the simple
version that locates bases at points in a control
space $K$ over $B$, and morphisms without incorporating
paths.

Let $K$ be a space,  $p \colon K \to B$ a 
(continuous) map
to a finite--dimensional
compact metric ANR space. We assume this
already here since the surgery sequence of (\cite{PQR})
requires these properties. 
Let $d \colon B \times B \to \mathbb{R}$ be a metric.
Let $\Lambda$ be a ring with involution and $1 \in \Lambda$,
i.e. $\Lambda$ is
a group ring.

A geometric module over $K$ is a free $\Lambda$--module
$M = \Lambda [S]$, $S$ a basis together with a map
$\varphi \colon S \to K$. It is required that for any $x \in K$,
$\varphi^{-1} (x) \subset S$ is finite. A morphism 
$f \colon M = \Lambda [S] \to N = \Lambda [T]$
is a collection $f_{st} \colon M_s \to N_t$
where $M_s = \Lambda [s]$, $N_t = \Lambda [t]$ such that for
a fixed $s \in S$ only finitely many $f_{st} \ne 0$.
The dual is
$M^* = \Lambda^* [S]$, where $\Lambda^* = \operatorname{Hom}_{\Lambda} (\Lambda, \Lambda)$,
so it
is essentially the same. However, if
$f \colon M \to N$ is a geometric morphism, its dual is
a geometric morphism $f^* \colon N^* \to M^*$.

Composition of geometric modules is defined
in the obvious way.
We define the radius of $f$ by 
$\operatorname{rad} (f) = \max \{ d(p\varphi (s), p\psi (t)) \; | \; f_{st} \ne 0 \}$. Here
$\psi \colon T \to K$ belongs to the geometric module $N$.
The map
$f$ is an $\varepsilon$--morphism if $\operatorname{rad} (f) < \varepsilon$.
The composition of an $\varepsilon$-- and a $\delta$--morphism
is an $(\varepsilon + \delta)$--morphism. The sum of an $\varepsilon$-- and  a 
$\delta$--morphism is a $\max \{ \varepsilon, \delta \}$--morphism.
See (\cite{RaYa1}) for further properties.

Chain complexes of geometric modules
are then defined as pairs $(C_{\mathop\#}, \partial_{\mathop\#})$ where
$\partial_n$ are geometric
morphisms. $(C_{\mathop\#}, \partial_{\mathop\#})$ is an
$\varepsilon$--chain complex if all $\partial_n$ are $\varepsilon$--morphisms.
A chain map $f \colon (C_{\mathop\#}, \partial_{\mathop\#}) 
\to (C_{\mathop\#}', \partial_{\mathop\#}')$
of geometric chain complexes is a $\delta$--chain
equivalence if there is a $\delta$--chain map
$g \colon (C_{\mathop\#}', \partial_{\mathop\#}') 
\to (C_{\mathop\#}, \partial_{\mathop\#})$ and chain homotopies
$\{ h_n \colon C_n \to C_{n+1} \}$, $\{ h'_n \colon C'_n \to C'_{n+1} \}$
of $g \circ f$ and $f \circ g$ with $\operatorname{rad} h_n < \delta$ and
$\operatorname{rad} h'_n < \delta$.
The composition of a $\delta$--chain equivalence with
a $\delta'$--chain equivalence is a $(\delta + \delta')$--chain
equivalence.
We observe that the dual $f^*$ of $f$ has the
same radius.

A chain equivalence is $\varepsilon$--contractible if
it is $\varepsilon$--chain equivalent to the zero
complex. If $f$ is an $\varepsilon$--chain equivalence
then its mapping cone is $3\varepsilon$--contractible.
An $\varepsilon$--mapping $f$ is a $2\varepsilon$--equivalence
if the mapping cone is $\varepsilon$--contractible.
A pair
$(C_{\mathop\#}, \psi_{\mathop\#})$ is called an $\varepsilon$--quadratic
geometric $\Lambda$--module complex if the maps
$\psi_s \colon C^{n-r-s} \to C_r$
have radius $< \varepsilon$. 
A pair
$(C_{\mathop\#}, \psi_{\mathop\#})$ is $\varepsilon$--Poincar\'e if
the mapping cone $\psi_0 + \psi_0^{\mathop\#}$ is $4\varepsilon$--contractible.
An $(n+1)$--dimensional $\varepsilon$--quadratic pair
$(C_{\mathop\#} \to D_{\mathop\#}, \psi_{\mathop\#}, \delta\psi_{\mathop\#})$ is a quadratic pair such
that
$
\delta\psi_s \colon D^{n+1-r-s} \to D_r
$
and
$
\psi_s \colon C^{n-r-s} \to C_r
$
have radius $< \varepsilon$. 
It is $\varepsilon$--Poincar\'e
if $C(f)^{n+1-r} = D^{n+1-r} \oplus C^{n-r} \to D_r$,
given by $(\delta\psi_0 + \delta\psi_0^{\mathop\#}, f \circ (\psi_0 + \psi_0^{\mathop\#}))$ has
$4\varepsilon$--contractible algebraic mapping cone
(note that it is a $2\varepsilon$--chain map).

Let $\delta \ge \varepsilon > 0$. Then
$L_n (p \colon K \to B, \varepsilon, \delta)$ is the set of
equivalence classes of
$n$--dimensional $\varepsilon$--quadratic $\varepsilon$--Poincar\'e
$\Lambda$--chain complexes on $p \colon K \to B$.
The equivalence relation is generated by
$\delta$--bordism defined in the obvious way.
It is actually shown that $\delta$--bordism is an equivalence relation.
The set $L_n (p \colon K \to B, \varepsilon, \delta)$ has a natural abelian
group structure given by direct sums.
As defined above,
these $\varepsilon$--$\delta$--$L$--groups
seem to not be calculable.

However, there are deep results identifying
these groups with homology groups in
certain spectra 
(\cite{PQR})
(\cite{PedYa})
(\cite{RaYa2}).
These spectra were constructed by Quinn
in (\cite{Qu0}) (see also (\cite{Ya})).
Here is a special case, 
due to Pedersen--Yamasaki (\cite{PedYa})
and Ranicki--Yamasaki (\cite{RaYa2}),
which is
most useful in $4$--dimensional surgery:

\begin{thm}
Consider $\Lambda = \mathbb{Z}$. Suppose
that
$p \colon K \to B$ 
is a fibration with simply connected
fibers. 
Then $L_n (p \colon K \to B, \varepsilon, \delta) \cong H_n (B, \mathbb{L})$,
for sufficiently small $\varepsilon$
and
$\delta$.
\end{thm}

Recall that $\mathbb{L}$ is the $4$--periodic (nonconnected)
surgery spectrum with $\pi_n (\mathbb{L}) = L_n (\{ 1 \})$.
A particular case is $p = \operatorname{Id} \colon B \to B$,
proved by
Pedersen--Quinn--Ranicki (\cite{PQR}) (a different proof was given by Ferry (\cite{Fe})):

\begin{thm}
Let $B$ be as above, then there is an
(assembly) isomorphism $L_n (B, \varepsilon, \delta) \cong H_n (B, \mathbb{L})$.
\end{thm}

We now come to the controlled surgery sequence
of (\cite{PQR}). First, we will explain the
surgery obstruction map.
Suppose $(f,b) \colon M^n \to X^n$ is a surgery problem,
$n = 2k$. Let $p \colon X \to B$ be a control map.
Here $X^n$ is an $n$--manifold or a $\delta$--Poincar\'e $n$--complex
over $B$ for sufficiently
small $\delta$, i.e.
$$
-\cap [X] \colon C^{*-n} (X) \to C_* (X)
$$
is a $\delta$--chain equivalence of $\Lambda$--modules over $B$, and the
cells of $X$ have diameter less than $\delta$ in $B$.
This holds for instance for generalized manifolds.

We want to describe the controlled surgery
obstruction of $(f,b)$ in $L_n (B, \varepsilon, \delta)$.
Let us assume that $p \colon X \to B$ is also $UV^1$.
The easiest way is to
consider $(f,b)$ as
an
element of
$H_n (X, \mathbb{L})$, and then its
image by $p_* \colon H_n (X, \mathbb{L}) \to H_n (B, \mathbb{L})$,
under the identification with $L_n (B, \varepsilon, \delta)$,
gives the controlled surgery obstruction.
It is however useful to write
down a controlled Wall obstruction,
i.e. in terms of a controlled $[K_k(f), \lambda,\mu]$.
However,
this must be an $\varepsilon$--quadratic
$\mathbb{Z}$--chain complex over $B$.

To obtain this, one does surgeries according
to the cell--structure of the relative complex
$(X,M)$, i.e. one substitutes $X$ by the
mapping cylinder of $f$.
In the first step one gets a normal cobordism of
$(f,b)$ to $(f',b') \colon M' \to X$ which is $(\varepsilon, k)$--connected,
i.e. any commutative diagram of continuous maps
\begin{center}
\input{diagram1.pst}
\end{center}
where $(L, L_0)$ is a $CW$--pair with $\dim L \le k$,
has an $\varepsilon$--controlled extension $\bar\alpha \colon L \to M'$
with $\bar\alpha |_{L_0} = \alpha_0$, and there is a homotopy
$h \colon L \times I \to X$ between $\alpha$ and $f' \circ \bar\alpha$,
with $\operatorname{radius} \{ ph(x,t) \; | \; t \in I \} < \varepsilon$ for
each $x \in L$.

We denote $(f',b') \colon M' \to X$ again by $(f,b) \colon M \to X$.
Note that $(\varepsilon, 2)$--connectedness is the  $UV^1 (\varepsilon)$--property.
Let $K_{\#}(f)$ be the kernel chain complex of $f$.
By the above we can assume that it is
$\varepsilon'$--chain equivalent to a geometric complex $E_{\#}$
over $B$ with $E_i = 0$ for $i \le k-1$, for
some $\varepsilon'$ depending on $\varepsilon$. We emphasize that
$K_{\#}(f)$ is the kernel complex of $M \rightarrow X$, not
of the universal covering.

Since $X$ is a $\delta$--Poincar\'e complex over $B$, $K_{\#}(f)$
has the structure of an $\varepsilon''$--quadratic geometric Poincar\'e
chain complex over $B$.
The next step is to apply controlled cell--trading 
and folding to get a chain complex $F_{\#}$
which is $\varepsilon'''$--equivalent to $E_{\#}$ and
$F_l = 0$ for $l \ne k$. For doing this,
one
needs that $p: X \to B$ is $UV^1$.
$F_{\#}$ is an $\varepsilon''''$--quadratic geometric Poincar\'e complex
with quadratic structure given by
intersection -- and self intersection numbers
$\lambda_{\mathbb{Z}}$, $\mu_{\mathbb{Z}}$ induced from $(f,b) \colon M \to X$.

The triple $(F_{\#}, \lambda_{\mathbb{Z}}, \mu_{\mathbb{Z}})$ represents
the controlled surgery obstruction of $(f,b)$.
If it is zero in $L_n (B, \varepsilon, \delta)$ then
controlled surgery can be completed if
$n \ge 5$ using the $UV^1$--property of $p \colon X \to B$
to find small Whitney  disks to
remove self--intersection numbers of
immersed spheres $S^k \to M^{2k}$,
representing generators in $F_k$. Completing
these surgeries one applies the controlled
Hurewicz--Whitehead theorem (\cite{Qu2}) (\cite{Qu0})
to get a controlled homotopy equivalence
$F'' \colon M'' \to X$.

This also works  in dimension $2k = 4$, since
one can apply the Controlled Disk Embedding
Theorem of Quinn ((\cite{Qu3}), cf. Disk Deployment
Lemma 3.2).
The following is the full statement from Pedersen--Quinn--Ranicki (\cite{PQR}):

\begin{thm}
Suppose that $B$ is as above. Then there is $\varepsilon_0 > 0$
such that for any $\varepsilon_0 > \varepsilon > 0$ there is $\delta > 0$
with the following property: If $X^n$ is a $\delta$--Poincar\'e
complex with respect to a $UV^1(\delta)$ map
$p \colon X^n \to B$ and $n \ge 4$, then there
is a controlled surgery exact sequence
$$
H_{n+1} (B, \mathbb{L}) \to \mathcal{S}_{\varepsilon, \delta} (X^n) \to
[X^n, G / _{TOP}] \mathop\to\limits^{\Theta} H_n (B, \mathbb{L}).
$$
\end{thm}

Here we must additionaly assume, that there
is a $TOP$
reduction of $\nu_X$.
Recall that
$\mathcal S_{\varepsilon, \delta} (X)$ consists of pairs $(M, f)$ where
$M^n$ is an $n$--manifold, $f \colon M \to X$ a $\delta$--homotopy 
equivalence over $p \colon X \to B$, modulo 
the equivalence relation: $(M, f) \sim (M', f')$ if
there is a homeomorphism $h \colon M \to M'$
such that $f$ and $f' \circ h$ are $\varepsilon$--homotopic
over $B$.
To define the Wall realization map
$H_{n+1} (B, \Bbb L) \to \Bbb S_{\varepsilon, \delta} (X)$,
one needs $\Bbb S_{\varepsilon, \delta} (X) \ne \emptyset$.

\bigskip\noindent
{\it Remark.}
Ranicki and Yamasaki worked
out, in a conceptual way, the controlled
surgery obstruction, using a controlled
version of the quadratic construction
(\cite{RaYa3}).

\bigskip\noindent
{\it Summary.}
Consider a surgery problem $(f,b) \colon M \to X$ with
control map $p \colon X \to B$. If $X$ is a $\delta$--Poincar\'e
complex for sufficiently small $\delta > 0$ over $B$, then
one can construct as above the controlled
surgery obstruction belonging to $L_n (B, \varepsilon, \delta)$.
If $p$ is additionally $UV^1 (\delta)$ for sufficientlly
small $\delta > 0$, then the controlled surgery sequence
holds.

\section{Some conclusions and comments}

In this section we present Quinn's approach and then
we consider Poincar\'e $4$--complexes with free fundamental groups.
We mentioned in the introduction Ranicki's
main result in high--dimensional surgery
theory: If $X^n$ is a Poincar\'e $n$--complex with
vanishing total surgery obstruction $s(X^n)$,
then $X^n$ is (simple) homotopy equivalent to
topological $n$--manifold $M^n$. Here $n \ge 5$.
One of the main objectives is to extend
this result
to
dimension $4$.

Here are  the key ideas of Quinn's approach
(\cite{Qu5}):
Let $X^4$ be a $4$--dimensional Poincar\'e complex.
\begin{itemize}
\item[(1)]
We investigate
the algebraic surgery sequence explained in $\S$2.
$$
\dots \to L_4(\pi_1(X^4))\to \mathbb{S}_4 (X^4) \to H_3 (X^4, \dot{\mathbb{L}}) \to
\dots $$
with $s(X^4) \in \mathbb{S}_4 (X^4).$
\item[(2)]
We consider the image of $s(X^4)$ under
the composite map
$$
\mathbb{S}_4 (X^4) \to H_3 (X^4, \dot{\mathbb{L}}) \to H_3 (X^4, \mathbb{L})
$$
and use the identification ($\S$3)
$$
H_3 (X^4, \mathbb{L}) \cong L_3 (X^4, \varepsilon, \delta), \forall \varepsilon < \varepsilon_0.$$
Thus $s(X^4)$ determines an element $[s(X^4)] \in L_3(X^4, \varepsilon, \delta)$,
i.e. 
$(D_{\mathop\#}, \psi_{\mathop\#}) = (S^{-1}C (-\cap[X^4]), \psi_{\mathop\#})$
as described
in $\S$2, carries an $\varepsilon$--quadratic
Poincar\'e structure, unique up to $\delta$--bordism.
\item[(3)]
If $[s(X^4)] = 0$, there is a $\delta$--null--bordism 
$(D_{\mathop\#}, \psi_{\mathop\#}) \to (E_{\mathop\#}, \delta \psi_{\mathop\#})$.
In fact, since we have assumed $s(X^4) = 0$,
$E_{\mathop\#}$ is contractible.
\item[(4)]
This bordism can be topologically
realized by a $\delta'$--homotopy equivalence
$X'^4 \to X^4$, where $X'^4$ is an $\varepsilon'$--Poincar\'e
$4$--complex. Here,
$(\delta', \varepsilon')$ depends on $(\delta, \varepsilon)$,
and becomes arbitrary small as $(\delta, \varepsilon)$ becomes
small. Ideas of surgery on Poincar\'{e} and normal spaces are 
used here (see
(\cite{Qu7})).
\item[(5)]
Choose a sequence $\{ \varepsilon_n \} \to 0$ 
and iterate the above construction to produce
a sequence $\{ X'^4_n \to X'^4_{n-1} \}_n$. Its limit
in the sense of (\cite{BFMW}) is an ANR  homology
4--manifold $X'^4$ which is homotopy equivalent
to $X^4$.
\end{itemize}

This approach can be summarized as follows:
Suppose $X^4$ is a Poincar\'e $4$--complex with
$s(X^4) = 0$. Then $X^4$ is (simple) homotopy
equivalent to an ANR  homology 4--manifold
$X'^4$.

\bigskip\noindent
{\it Remarks.} 
(1)
The topological
realization step ($4$) requires a highly
$\delta$--connected null bordism $(E_{\mathop\#}, \delta \psi_{\mathop\#})$, which
is {\sl not} guaranteed
when $n$ is even.

(2)
Starting with a relative Poincar\'e
complex $(X^4, \partial X^4)$ such that $\partial X^4$ is
a topological 4--manifold, the $4$--dimensional
resolution theorem (\cite{Qu3}) implies that
$X'^4$ is a topological 4--manifold.

\bigskip
For the rest of this section we consider
Poincar\'e
$4$--complexes $X^4$ with free nonabelian
fundamental groups, i.e. $\pi_1(X^4) \cong \mathop{*}\limits_1^p \mathbb{Z}$.
We benefit
from
the special topology of such
complexes, in particular:

\begin{thm}
\begin{itemize}
\item[(a)]
$X^4$ is (simple) homotopy equivalent to
$\{ \mathop\vee\limits_1^p (S^1 \vee S^3) \vee (\mathop\vee\limits_1^q S^2) \} 
\mathop\cup\limits_{\varphi} D^4$; 
and
\item[(b)]
If the $\Lambda$--intersection form
$$
\lambda_{\Lambda} \colon H_2 (X^4, \Lambda) \times H_2 (X^4, \Lambda) \to \Lambda
$$
is extended from the $\mathbb{Z}$--intersection form
$$
\lambda_{\mathbb{Z}} \colon H_2 (X^4, \mathbb{Z}) \times H_2 (X^4, \mathbb{Z}) \to \mathbb{Z}
$$
then $X^4$ is (simple) homotopy equivalent to
$Q^4 \# M'^4$, where $Q^4 = \mathop\#\limits_1^p (S^1 \times S^3)$, and $M'^4$ is
a simply connected topological 4--manifold
determined by $\lambda_{\mathbb{Z}}$.
\end{itemize}
\end{thm}

For proofs see 
(\cite{HegPic})
(\cite{HegRepSpagg})
(\cite{Hil})
(\cite{MatKat}). 
We note here that the first Postnikov invariant for $X^4$ vanishes.
Theorem 4.1 implies (what is much easier
to see) that there is a degree 1 map
$
p \colon X^4 \to Q^4.
$

\begin{lem}
The assembly maps satisfy the following properties:
\begin{itemize}
\item[(a)]
$A \colon H_4 (X^4, \dot{\mathbb{L}}) \to L_4(\pi_1(X^4))$
is onto; and
\item[(b)]
$A \colon H_3 (X^4, \dot{\mathbb{L}}) \to L_3(\pi_1(X^4))$
is injective.
\end{itemize}
\end{lem}

\noindent
{\it Proof.}
Assembly is a natural
construction so we have the
commutative diagram
\begin{center}
\input{diagram2.pst}
\end{center}
A spectral sequence argument shows that
$p_* \colon H_4 (X^4, \dot{\mathbb{L}}) \to H_4 (Q^4, \dot{\mathbb{L}})$ is onto,
and $p_* \colon H_3 (X^4, \dot{\mathbb{L}}) \to H_3 (Q^4, \dot{\mathbb{L}})$ is an isomorphism.
If $B = B(\pi_1 (X^4))$ is the classifying space,
$c \colon Q^4 \to B$ the classifying map,
then
$c_* \colon H_l (Q^4, \dot{\mathbb{L}}) \to H_l (B, \dot{\mathbb{L}})$ is an
isomorphism for $l = 3,4$ by similar
arguments. However, for free fundamental
groups,
$A \colon H_l (B, \dot{\mathbb{L}}) \to L_l (\pi_1(x))$
is an isomorphism. This proves
the lemma.

\begin{cor}
If $X^4$ is a Poincar\'e $4$--complex,
then $s(X^4)$ is zero. In fact, the same holds for the
algebraic structure set $\mathbb{S}_4 (X^4) = \{0\}$.
\end{cor}

\noindent
{\it Proof.}
This follows from Lemma 4.2
and  the algebraic surgery sequence
$$
\to H_4 (X^4, \dot{\mathbb{L}}) \mathop\to\limits^A L_4(\pi_1(X^4)) \to \mathbb{S}_4(X^4) \to
$$
$$
\to H_3(X^4,\dot{\mathbb{L}}) \mathop\to\limits^A L_3(\pi_1(X^4)) \to \cdots
$$

\bigskip
By the  discussion
above it is plausible to conjecture:

\begin{conj}
Any Poincar\'e $4$--complex $X^4$, such that
$\pi_1(X^4) \cong \mathop{*}\limits_1^p \mathbb{Z}$,
is (simple) homotopy equivalent to an
ANR  homology 4--manifold.
\end{conj}

Part (b) of the above theorem confirms
Conjecture 4.4 for the case when
$\lambda_{\Lambda}$ is extended
from $\lambda_{\mathbb{Z}}$. Indeed, in this case $X^4$ is
homotopically a manifold.
In general case we obtain a "stable" result:

\begin{cor}
If $X^4$ has a
free nonabelian
fundamental group, then $X^4 \#
(\mathop\#\limits_1^r S^2 \times S^2)$
is (simple) homotopy equivalent to a topological
4--manifold.
\end{cor}

\noindent
{\it Proof.}
Since $s(X^4) = 0$, there is
a degree 1
normal map $(f,b) \colon M^4 \to X^4$
whose
Wall obstruction is
zero.
This means that $(K_2(f), \lambda, \mu)$
is stably hyperbolic.
The result then follows from (\cite{HegRepSpagg}).

\bigskip\noindent
{\it Remark.}
$X^4 \# (\mathop\#\limits_1^r S^2 \times S^2)$ is the connected
sum
made inside a $4$--cell in
$X^4$.
\bigskip

The controlled surgery method also works for certain other fundamental groups. 
We  have proved this for those Poincar\'{e} complexes whose 
fundamental 
group is that of a torus knot (\cite{HegRep}):

\begin{thm}
Let $X^4$ be a 4-dimensional Poincar\'{e} complex such that
$\pi_{1}(X^4) \cong \pi_1(S^3\setminus K)$, where $K \subset S^3$ is a torus knot, and suppose that 
 $s(X^4)=0$. Then $X^4$ is (simple) homotopy equivalent to a closed topological
$4$--manifold. 
\end{thm}

Yamasaki (\cite{Ya2})
has recently
proved that Theorem 4.6 holds also for hyperbolic knots $K \subset S^3$. Note that in order to verify Theorem 4.6 for {\sl all} knots  $K \subset S^3$ it would suffice, by Thurston's theorem, to answer in affirmative the following question:

\begin{quest}
Does Theorem 4.6 hold also if  $K \subset S^3$ is a satellite knot? 
\end{quest}


\end{document}